\documentclass[a4paper,11pt]{article}

\usepackage{cite}

\usepackage{amsmath, amsthm, amsfonts, amssymb, dsfont, color}
\usepackage{srcltx}

\usepackage{graphicx}

\theoremstyle{plain}
\newtheorem{theorem}{Theorem}[section]

\newtheorem{lemma}[theorem]{Lemma} 

\newtheorem{corollary}[theorem]{Corollary}

\theoremstyle{remark}

\theoremstyle{definition}

\numberwithin{equation}{section}

\def\e{\varepsilon}

\def\G{\Gamma}

\def\p{\partial}
\def\D{\Delta}

\def\E{\mbox{\rm e}}
\def\a{\alpha}
\def\b{\beta}

\def\z{\zeta}

\def\Odr{\mathcal{O}}
\def\H{W_2}

\def\iu{\mathrm{i}}

\DeclareMathOperator{\RE}{Re}

\DeclareMathOperator{\spec}{\sigma}

\DeclareMathOperator{\dist}{dist}

\newcommand{\RR}{\mathds{R}}
\newcommand{\ZZ}{\mathds{Z}}
\newcommand{\NN}{\mathds{N}}
\newcommand{\CC}{\mathds{C}}

\renewcommand{\varpi}{\omega}

\begin{document}

\allowdisplaybreaks

\begin{center}
\textbf{\Large Gap opening and split band edges in waveguides coupled
by a periodic system of small windows}

\medskip

{\large D. I. Borisov\footnote{D. Borisov was partially supported by RFBR grant no. 10-01-00118 and Federal Task Program ``Scientific and pedagogical staff of innovative Russia for 2009-2013 гг.'', contract no. 02.740.11.0612, and by FCT grant no. ptdc/mat/101007/2008.} (M. Akmulla Bashkir State Pedagogical University and Institute of mathematics of the Russian Academy of Sciences, Ufa, Russia,\\ E-mail: BorisovDI@yandex.ru)\\[\medskipamount]
K. V. Pankrashkin (Unit\'e mixte de recherche 8628, CNRS and University Paris-Sud 11, Orsay, France, E-mail: konstantin.pankrashkin@math.u-psud.fr)}

\end{center}

\begin{abstract}

At the example of two coupled waveguides we construct a periodic second order differential operator
acting in a Euclidean domain and having spectral gaps whose edges are attained
 strictly inside the Brillouin zone. The waveguides are modeled by the Laplacian in two infinite strips
of different width that have a common interior boundary. On this common boundary we impose the Neumann boundary condition
but cut out a periodic system of small holes, while on the remaining exterior boundary we impose the Dirichlet boundary condition.
It is shown that, by varying the widths of the strips and the distance between the holes,
one can control the location of the extrema of the band functions as well as the number of the open gaps.
We calculate the leading terms in the asymptotics for the gap lengths and the location of the extrema.

\medskip

\end{abstract}

\section*{Introduction}

It is well known that a large class of periodic operators have a band spectrum, i.e.
the spectrum is the union of finite segments. One of the important questions in the study of such operators
concerns the existence of spectral gaps, i.e. open intervals lying outside the spectrum but whose
edges are in the spectrum. The study of this property is motivated by various applications, in particular,
by the problems arising in the theory of photonic crystals \cite{PK}.

The one-dimensional Schr\"odinger operator with a periodic potential (Hill operator)
is a classical example of a periodic operator. The only possible gap edges for these operators
are the values of the energy for which the eigenvalue equation has a non-trivial periodic or anti-periodic solution.
While this reflects some specific features of the operator, the same was implicitly assumed to be true
for more general operators; in those cases the periodic and anti-periodic eigenvalues
are replaced by the values of the band functions at some special values of the quasimomentum (see below).
A discussion of analytic aspects of this property appeared just quite recently
in the works \cite{EKW,HKSW}, which served the motivation for the present paper.
The work \cite{HKSW} contains, in particular, several examples of quite sophisticated multi-dimensional
combinatorial operators whose gap edges are attained at intermediate values of the quasimomentum,
and the work \cite{EKW} is devoted to a more detailed analysis of $\ZZ$-periodic operators.

It should be emphasized that similar effects, like attaining extremal values at intermediate values of quasimomenta
and the presence of multiple extrema, are known in the physics literature for a long time as split band edge,
and these effects play an essential role, in particular, in slowing waves in photonic crystals,
see e.g. \cite{FV,PH1,PH2}. At the physics level of rigor it was noted that
the split band edges can appear when coupling several $\ZZ$-periodic structures with different
characteristics, cf. \cite{PH2}. In the present work we give this assertion a rigorous sense
at the example of two non-symmetric waveguides. We consider the systems consisting of
two two-dimensional strips of different width coupled by a periodic system
of windows. The Dirichlet boundary conditions are imposed at the exterior boundary,
and the Neumann boundary conditions at the interior one.
It is shown that, if some inequalities between the period and the waveguide widths hold,
under the assumption that the windows are small, the associated operator (Laplacian) has gaps
whose edges are attained at intermediate values of quasimomentum.
We note that we use explicitly the $\ZZ$-periodicity (i.e. periodicity in one direction
and compactness in the others), and we do not expect the results to be valid for
systems with multi-dimensional periodicity (see e.g. the review \cite{HP}
for possible ways of gap opening in multidimensional systems).

The questions of gap opening in coupled and perturbed waveguides (in situations which are different from the one under consideration)
were already studied in numerous works. In particular, the existence of gaps
in the spectrum of a periodically curved waveguide was shown in \cite{Yo2},
and the work \cite{FS} deals with the study of the asymptotics of the bands and the gaps
of curved waveguides with respect to the cross-section size.
Several works studied by various methods the waveguides represented as a chain of domains
coupled by windows or links \cite{sn2,sn3,KP,Yo}; in that case the number of gap
grows to infinity. The works \cite{sn,CNP,naz2} are the closest to us from the point of view of methods,
they show the gap opening in cylindrical waveguides with periodic
localized perturbation at the boundary. Furthermore, there were some works in the physics literature \cite{pop1,pop2,pop3},
dealing with the study of gaps in waveguides coupled or perturbed by periodically varying boundary conditions.
The main contribution of the present work consists in a rigorous description
of the split band edge effect in coupled waveguides, which, to our knowledge, was not found earlier.

The paper is organized as follows. In section \ref{sec1} the rigorous problem setting is given
and the most important result, Theorem \ref{th1.1}, is formulated, whose proof is given in the two subsequent sections.
In section \ref{sec2} we are concerned with an upper bound of the shift of band functions of uncoupled and coupled waveguides.
In section \ref{sec3} the crossings of the unperturbed band functions are analyzed in greater details,
and the asymptotics of the perturbed band functions near these crossings is given.

\section{Problem setting}\label{sec1}

For the sake of convenience let us denote $d_+:=\pi$, $d_-=d>0$, and, furthermore,
let $h$ and $\e$ be parameters satisfying
$0\leqslant \e<h$. In $\RR^2$ we consider infinite straight strips
\begin{equation*}
\Pi_+:=\RR\times(0,d_+),\quad
\Pi_-:=\RR\times(-d_-,0),\quad
\Pi:=\RR\times(-d_-,d_+),
\end{equation*}
and introduce the operator $H^\e$ acting in $L_2(\Pi)$
as the Laplacian with the Dirichlet boundary conditions at the lines $x_2=\pm d_\pm$ and the Neumann boundary conditions
at the line $x_2=0$ except at the intervals $\varpi^\e_n:=(2nh-\e,2nh+\e)\times\{0\}$, $n\in\ZZ$.
It is clear that for $\e=0$ we just have the equality $H^0=H_+\oplus H_-$, where
$H_\pm$ are the Laplacians in the waveguides $\Pi_\pm$ with the Dirichlet boundary conditions at the lines
$x_2=\pm d_\pm$ and the Neumann boundary conditions at the line $x_2=0$.

\begin{figure}
\centering

\includegraphics[width=100mm]{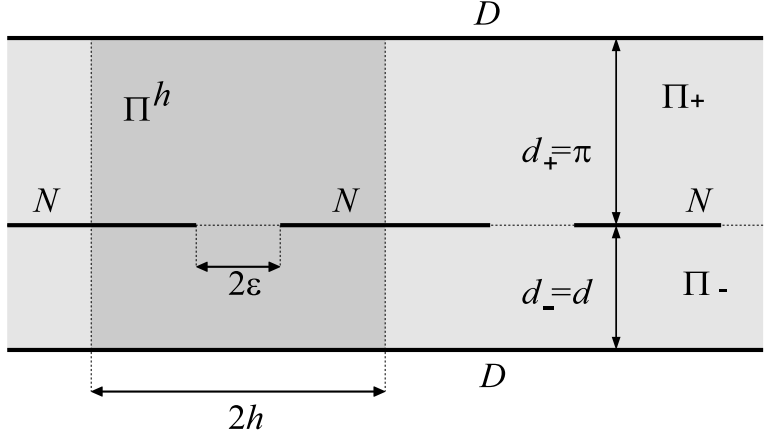}

\caption{The waveguides $\Pi_+=\RR\times(0,\pi)$ and $\Pi_-=\RR\times(-d,0)$ coupled by a periodic system of windows.
The parts of the boundary on which the Dirichlet and Neumann boundary conditions
are imposed are marked by the symbols $D$ and $N$ respectively.
The elementary cell $\Pi^h$ is highlighted. We denote $\Pi_\pm^h:=\Pi^h\cap \Pi_\pm$.}\label{fig1}
\end{figure}

Therefore, the operator $H^0$ is the Hamiltonian of the system consisting of
two non-interacting waveguides $\Pi_+$ and $\Pi_-$, while
$H^\e$ describes the same waveguides coupled by the
$2h$-periodic system of windows $\varpi^\e_n$ of length $2\e$, see Figure~\ref{fig1}.

The present works analyses the situation of non-symmetric waveguides: $d\ne \pi$.
Without loss of generality we assume
\begin{equation}\label{1.1}
d<\pi.
\end{equation}
The spectra of the operators $H_\pm$ are well studied. There holds
\begin{equation*}
\spec (H_\pm)=\left[\dfrac{\pi^2}{4d_\pm^2},+\infty\right),
\end{equation*}
hence, under the assumption (\ref{1.1}), the spectrum of  $H^0$ is the half-axis $[1/4,+\infty)$.
Moreover, the spectra of $H_\pm$ and $H^0$ are absolutely continuous.

Let us denote by $\Pi^h$ the elementary cell given by $\Pi^h:=\Pi\cap\{x=(x_1,x_2): |x_1|<h\}$.
The periodicity with respect to $x_1$ and the application of the Bloch-Floquet theory
allows one to decompose $H^\e$ into a direct integral:
\begin{equation}\label{1.2}
H^\e\simeq\int^\oplus_{(-\pi,\pi]} H^\e(k) dk,
\end{equation}
where $H^\e(k)$ is the operator acting in $L_2(\Pi^h)$ by
\begin{equation}
       \label{eq-hk}
H^\e(k)=\Big(\iu\dfrac{\partial}{\partial x_1}-\dfrac{k}{2h}\Big)^2 - \dfrac{\partial^2}{\partial x_2^2},
\end{equation}
with the boundary conditions
\begin{subequations}
  \label{1.10}
\begin{gather}
 \label{1.10a}
u=0\quad\text{on}\quad \G_\mathrm{D}^h,
\\
  \label{1.10b}
\frac{\p u}{\p x_2}=0\quad \text{on}\quad \G_\mathrm{N}^{h,\e},\\
  \label{1.10c}
u|_{x_1=-h}=u|_{x_1=h},\quad \frac{\p u}{\p x_1}\Big|_{x_1=-h}=\frac{\p u}{\p x_1}\Big|_{x_1=h},
\end{gather}
where
\begin{align*}
\G_\mathrm{D}^h&:=\big\{x=(x_1,d_+): |x_1|<h\big\}\cup\big\{x=(x_1,-d_-): |x_1|<h\big\},\\
\G_\mathrm{N}^{h,\e}&:=\big\{x=(x_1,0): \e<|x_1|<h\big\}.
\end{align*}
\end{subequations}
The variable $k$ in the preceding formulas is usually referred to as the quasimomentum
or the Bloch parameter.

The spectrum of each of the operators $H^\e(k)$ is discrete. Let us denote by
$E^\e_l(k)$, $l\in\NN$, its eigenvalues taken with their multiplicities and enumerated
in the non-decreasing order,
$E^\e_l(k)\leqslant E^\e_{l+1}(k)$.
The functions $k\mapsto E^\e_l(k)$, called band functions
(or dispersion laws, or dispersion relations),
are $2\pi$-periodic, continuous and even (in virtue of the real-valuedness of the Laplacian),
and the spectrum of $H^\e$ is nothing but the union of the ranges of these functions,
\[
\spec(H^\e)=\bigcup_{l\in\NN} B_l^\e,
\quad
B_l^\e=E^\e_l\big((-\pi,\pi]\big).
\]
The segment $B_l^\e=:[\a_l^-,\a_l^+]$  is called $l$th band
of the operator $H^\e$. If for some $l\in\NN$ one has $\a_l^+<\a_{l+1}^-$,
then the interval $(\a_l^+,\a_{l+1}^-)$ is called a gap of the operator $H^\e$.

The central question of our work is the question on the existence of gaps in the spectrum of $H^\e$ for small $\e$,
as well as the question on the values of the quasimomentum $k$ at which the extremal values $\a_l^\pm$
of the dispersion laws are attained. In particular, the attaining the extremal values
outside the particular values $k=0$ and $k=\pm\pi$ will be put into evidence.

To formulate the main results we need some additional notation.
Let us denote $\Pi_\pm^h:=\Pi_\pm\cap\{x: |x_1|<h\}$. In the spaces $L_2(\Pi_\pm^h)$ consider the operators $H_\pm(k)$
given by the differential expression (\ref{eq-hk}) with the boundary conditions (\ref{1.10}), where in the Neumann condition (\ref{1.10b})
we formally assume $\G_\mathrm{N}^{h,0}:=\{x=(x_1,0): |x_1|<h\}$.
Analogously to (\ref{1.2}), one has the decomposition
\begin{equation*}
H_\pm\simeq\int^\oplus_{(-\pi,\pi]} H_\pm(k) dk.
\end{equation*}

The operators $H_\pm(k)$ have discrete spectra, and their eigenvalues and eigenfunctions
are calculated explicitly. The eigenvalues are
\begin{equation}
     \label{eq15}
E^\pm_{m,p}(k)=\Big(\dfrac{k+2\pi m}{2h}\Big)^2 + \Big(\dfrac{\pi}{d_\pm}\Big)^2\Big(p+\frac{1}{2}\Big)^2,
\quad m\in\ZZ,\quad p\in\ZZ_+\equiv\NN\cup\{0\},
\end{equation}
for $k\notin \{0,\pi\}$ the eigenvalue $E^\pm_{m,p}(k)$ is simple with the eigenfunction
\begin{equation}
 \label{eq-ef}
\psi_{m,p}(x,y)=\E^{\frac{\iu \pi m x}{h}} \cos \Big[\dfrac{\pi}{d_\pm}\Big(p+\frac{1}{2}\Big)y\Big],
\end{equation}
while for $k\in\{0,\pi\}$ the respective eigenvalue is double.
The spectrum of $H^0(k)$ is exactly the union, taking into account the multiplicity,
of the eigenvalues $E^+_{m,p}(k)$ and $E^-_{m,p}(k)$ for all $m\in\ZZ$, $p\in\ZZ_+$.

Let us note that the multiplicity grows infinitely at high energies: each point of the spectrum $E$
of the operator $H^0$ belongs the the ranges of $\kappa(E)$ function $E^\pm_{m,p}$, and
$\kappa(E)\xrightarrow{E\to+\infty} +\infty$.
Our main result is contained in the following theorem and describes the behavior
of the spectrum of $H^\e$ in the neighborhood of the points $E$ with $\kappa(E)\leqslant 4$.
Note that the inequality $\kappa(E)\leqslant 4$ implies automatically $E<\dfrac{9}{4}$.

\begin{figure}

\centering

\begin{tabular}{cc}
\includegraphics[height=100mm]{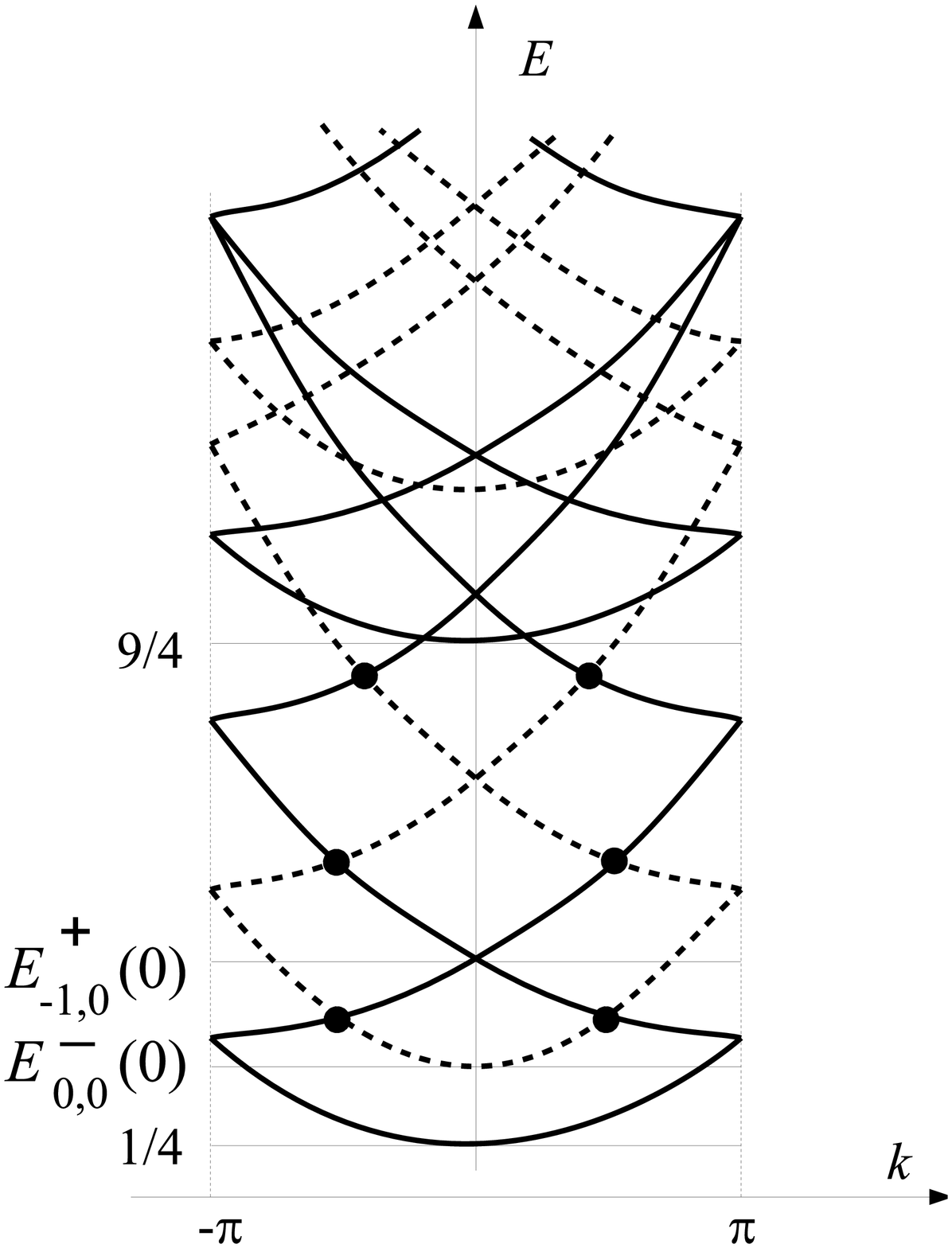}
&
\includegraphics[height=100mm]{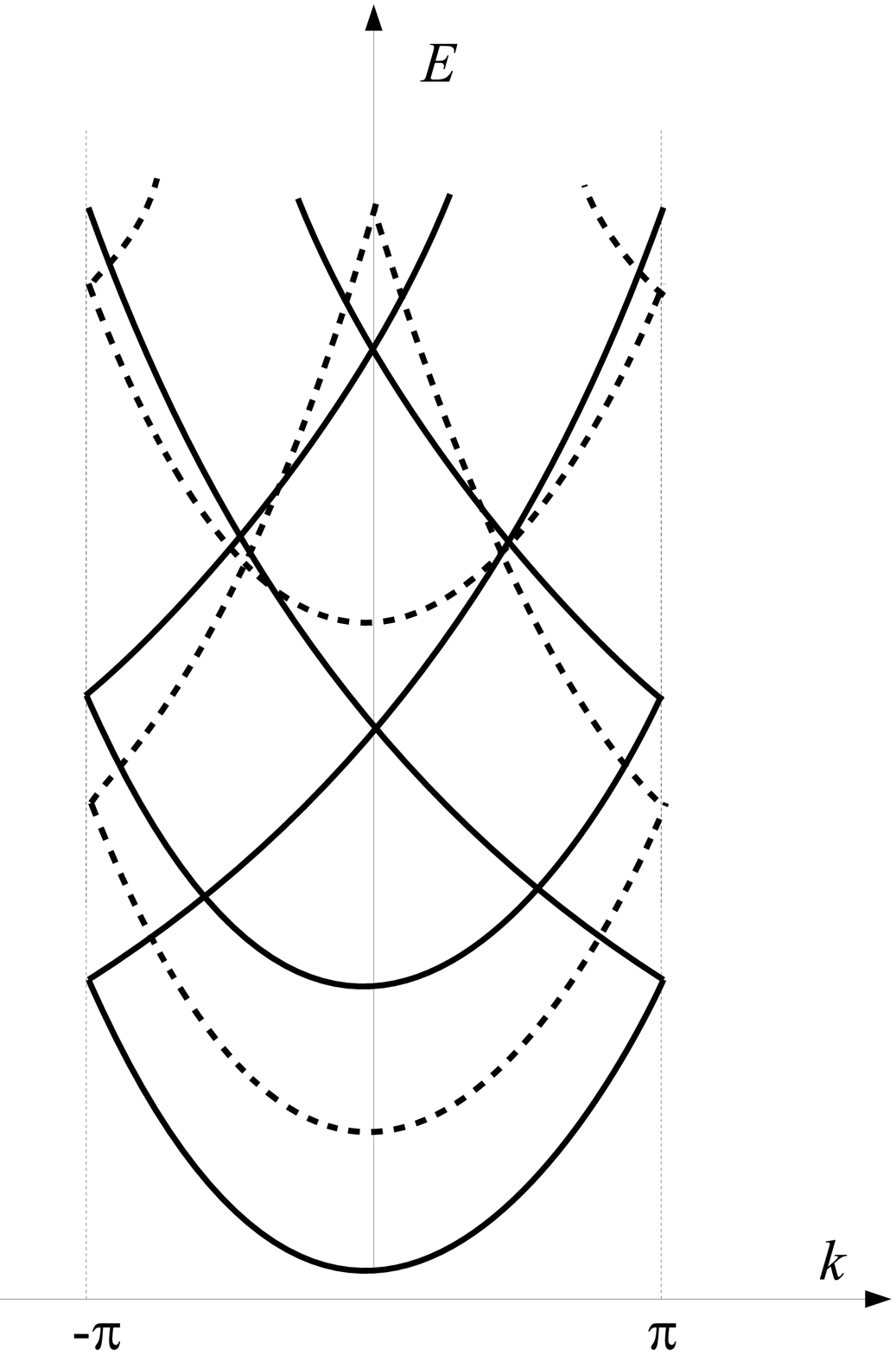}
\\
(a) & (b)
\end{tabular}

\caption{A schematic view of two possible configurations of band functions for $H^0$. The solid lines correspond to the functions
$E^+_{m,p}$, the dashed ones correspond to the functions $E^+_{m,p}$. (a) The intersections of the band functions satisfying
the assumption of Theorem \ref{th1.1} are marked by fat points. (b) There are no points satisfying the assumptions
of Theorem \ref{th1.1}.}\label{fig2}
\end{figure}

\begin{theorem}\label{th1.1}
Assume that the inequality (\ref{1.1}) holds and that there exists a point $k_0\in(0,\pi)$ and $m,n\in\ZZ$
for which one has
\begin{subequations}
\label{1.4}
\begin{gather}
       \label{1.4a}
E_0:=E^+_{n,0}(k_0)=E^-_{m,0}(k_0)<\dfrac{9}{4},\\
\label{1.4b}
\dfrac{\p E^+_{n,0}}{\p k}(k_0)\cdot\dfrac{\p E^-_{m,0}}{\p k}(k_0)<0.
\end{gather}
\end{subequations}
Then, for sufficiently small $\e$, there is a gap $\big(\a_l(\e), \a_r(\e)\big)$ in the spectrum
of the operator $H^\e$ whose endpoints have the asymptotics
\begin{equation}\label{1.3}
\a_j(\e)=E_0-\frac{\tau_j}{4h|\ln\e|}+\Odr(|\ln\e|^{-2}), \quad j=l,r,
\end{equation}
where
\begin{equation}
\label{1.6}
\begin{aligned}
&\tau_l:=2\sqrt{\z(1-\b^2)}-\b(\z-1)-\z-1,
\\
&\tau_r:=-2\sqrt{\z(1-\b^2)}-\b(\z-1)-\z-1,
\\
&\z:=\frac{\pi}{d}, \quad \b:=\frac{\pi(m+n)+k_0}{\pi(n-m)}.
\end{aligned}
\end{equation}
The band functions $E_l^\e(k)$ and $E_r^\e(k)$ describing the points of the spectrum to the left, respectively to the right,
from the gap $\big(\a_l(\e),\a_r(\e)\big)$, attain their extremal values
$\a_l(\e)\equiv \max_{|k|\leqslant\pi} E_l^\e(k)$ and $\a_r(\e)\equiv \min_{|k|\leqslant\pi} E_r^\e(k)$
at the points $\pm k_l(\e)$ and $\pm k_r(\e)$, respectively:
\begin{equation*}
E_l^\e\big(\pm k_l(\e)\big)=\a_l(\e),\quad E_r^\e\big(\pm k_r(\e)\big)=\a_r(\e).
\end{equation*}
One has the asymptotics
\begin{equation}
\label{1.7}
\begin{gathered}
k_{l/r}(\e)=k_0+ \frac{\sigma_{l/r}}{\ln\e}+
\Odr(|\ln\e|^{-3/2}),\\
\sigma_l=-\frac{\b h}{\pi(n-m)}
\sqrt{\frac{\z}{1-\b^2}}+\frac{(1-\z)h}{2\pi(n-m)},\\
\sigma_r=\frac{\b h}{\pi(n-m)}\sqrt{\frac{\z}{1-\b^2}}+\frac{(1-\z)h}{2\pi(n-m)}.
\end{gathered}
\end{equation}
\end{theorem}

Let us remark that one can use the methods of the works \cite[Ch. 9]{GN}, \cite{MP, G} and proceed with a deeper analysis
in order to obtain further terms of the asymptotics \eqref{1.7} as some powers of $1/\ln\e$.
Moreover, one can show some analyticity properties with respect to $1/\ln\e$,
see \cite{hc,G,GN}.

Note that the theorem stated does not guarantee that the inequalities (\ref{1.4})
hold for at least one combination $(m,n,k_0)$. Moreover, one can easily choose the parameters
in such a way that this condition will not be satisfied at any point;
for example, this can be obtained by taking $h$ sufficiently small
(see Fig.~\ref{fig2} for illustration of possible configurations).
On the other hand, to satisfy the conditions (\ref{1.4}) in at least one point it is sufficient
to have the inequality $E^-_{0,0}(0)<E^+_{-1,0}(0)<\dfrac{9}{4}$,
see Fig.~\ref{fig2}(b).
Rewriting this inequality with the explicit expressions (\ref{eq15}) we arrive at

\begin{corollary}
Let the conditions
\begin{equation*}
h>\dfrac{\pi}{\sqrt{2}}
\quad \text{and} \quad
\dfrac{\pi}{\sqrt{\Big(\dfrac{2\pi}{h}\Big)^2+1}} < d <\pi
\end{equation*}
hold. Then there exists $k_0\in(0,\pi)$ for which  (\ref{1.4}) holds
with $n=-1$ and $m=0$, and the operator $H^\e$,
for $\e$ sufficiently small, has a gap lying in an
$\Odr(\ln^{-1}\e)$-neighborhood of the respective point
$E_0$, and the adjacent band functions attain the respective extremal values
in an $\Odr(\ln^{-1}\e)$-neighborhood of  $k_0$.
\end{corollary}

\begin{figure}

\centering

\includegraphics[width=120mm]{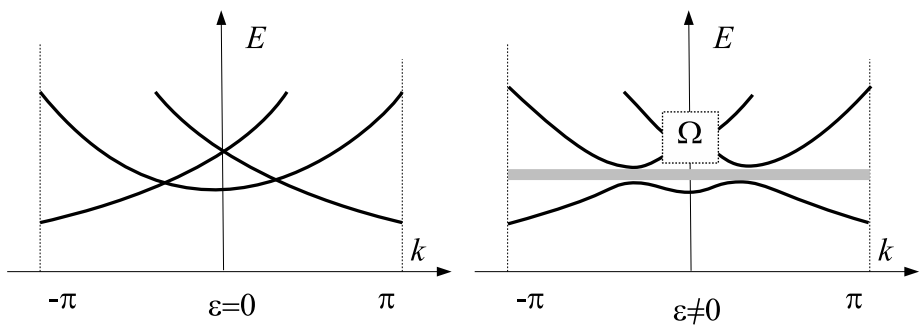}

\caption{Gap opening in the situation of Theorem \ref{th1.1}.}\label{fig3}

\end{figure}

Figure \ref{fig3} illustrates the behavior of the band functions near the points
at which the assumptions of Theorem \ref{th1.1}.
Note that even though the technique of proof of Theorem~\ref{th1.1} can be transferred without significant
additional efforts to the case $k\in\{0,\pi\}$ and allows one
to clarify the behavior of the band functions at such crossings (domain $\Omega$ on Fig.~\ref{fig3})
by some additional computations, these parts are overlapped by the other parts
of the graphs when projecting to the axis $E$, hence
these crossings do not participate in the gap opening.

We remark that by a suitable choice of the parameters $d$ and $h$ one can satisfy the assumptions of Theorem \ref{th1.1}
at any prescribed value of quasimomentum $k_0\in (0,\pi)$. In a similar way one can
obtain any prescribed number of crossings in the interval $(1/4, 9/4)$, and each of them will satisfy the conditions (\ref{1.4}).
Each of these crossings will generate its own gap by Theorem~\ref{th1.1}. Hence by fitting the parameters $d$ and $h$
one can open any prescribed number of gaps on the interval $(1/4, 9/4)$,
and each of them is covered by Theorem~\ref{th1.1}. 
The respective values  $k_0$ associated with different gaps will, in general, also be different.

\section{Estimating the shift of band functions}\label{sec2}

One of the main ingredients of the proof of Theorem~\ref{th1.1} comes from
a certain rather general estimate of the difference $\big|E^\e_l(k)-E^0_l(k)\big|$.
This is discussed in the present section. The estimate needed is contained in the following Lemma.

\begin{lemma}\label{prop1}
For any $l\in\NN$ there exists $C_l>0$ such that, for $\e>0$ sufficiently small,
and all $k\in(-\pi,\pi]$ one has the inequality
\begin{equation}
 \label{eq-ele}
0\leqslant E_l^\e(k)-E_l^0(k)\leqslant \dfrac{C_l}{|\ln \e|}.
\end{equation}
\end{lemma}

For the proof of the required estimates, like how it was done
in similar situations in \cite{sn2,sn3}, we use the max-min principle.
Let us consider the quadratic forms $a^0_k$ and $a^\e_k$ associated with the operators $H^0(k)$ and
$H^\e(k)$  respectively. The standard constructions show that the form
$a^0_k$ is given by the expression
\[
a^0_k(u,u)=\int_{\Pi^h}
\bigg(
\Big|
\Big(\iu\dfrac{\p}{\p x_1}-\dfrac{k}{2h}\Big) u(x_1,x_2)
\Big|^2
+
\Big|
\dfrac{\p u(x_1,x_2)}{\p x_2}
\Big|^2
\bigg)
dx_1\,dx_2
\]
on the domain $D^0$ consisting of the functions
$u\in W^1_2(\Pi_+^h)\oplus W^1_2(\Pi_-^h)$
satisfying the boundary conditions
\begin{equation}
    \label{eq-ae0}
u|_{\Gamma^h_D}=0,
\quad
u|_{x_1=-h}=u|_{x_1=h}.
\end{equation}
In turn, the form $a^\varepsilon_k$ is the restriction of $a_k^0$
to the set $D^\varepsilon$ consisting of the functions $u\in D^0$ satisfying,
in addition to (\ref{eq-ae0}), the matching condition on the window
\begin{equation}
   \label{eq-ae1}
u^+|_{\omega^\e_0}=u^-|_{\omega^\e_0},
\end{equation}
where $u^\pm$ stands for the restriction of $u$ onto $\Pi_\pm^h$.

According to the max-min principle \cite[Th. XIII.2]{RS4}, for the eigenvalues $E^\varepsilon_l(k)$ one has the equality
\begin{equation}
      \label{eq-minmax}
E^\e_l(k)=\max_{L\in S_{l-1}} \min_{\substack{u\in D^\e\cap L\\ u\ne 0}}\dfrac{a^0_k(u,u)}{\|u\|^2_{L^2(\Pi^h)}},
\quad l\in\NN, \quad \e\geqslant 0,
\end{equation}
where we denote by $S_m$ the set of all closed subspaces of $L^2(\Pi^h)$ of codimension $m$,
(i.e. $\dim L^\perp=m$ for all $L\in S_m$).
As  $D^\e\subset D^0$, the left-hand part of the inequality (\ref{eq-ele}) follows directly from \eqref{eq-minmax}.

To prove the right-hand part of (\ref{eq-ele}) 
we use the approach similar to the one used in \cite{sn2,sn3}.
Introduce first a family of cut-off functions.
Let $\psi:[0,+\infty)\to\RR$ be an infinitely smooth function with $\psi(s)=1$ for $s\leqslant \frac 1 2$
and $\psi(s)=0$ for $s\geqslant 1$. For $\e>0$ let us introduce the functions $\Phi_\e:\RR^2\to\RR$ by the equalities
\[
\Phi_\e(x)=\psi\Big(\Big|\dfrac{\ln |x|}{\ln\e}\Big|\Big),\quad x\ne 0, \qquad
\Phi_\e(0)=0.
\]
Note that for any function $v\in D^0$ the product $u=\Phi_\e v$ for $\e$ sufficiently small
is still in $D^0$ and, moreover, vanishes near the window $\omega^\e_0$. 
Therefore one has the equality \eqref{eq-ae1} and the inclusion $u\in D^\e$.
Moreover, denoting $B_1=\max_{s\in\RR}\big|\psi'(s)\big|$,
for $\e$ sufficiently small we have the estimates
\begin{equation}
   \label{eq-ph1}
\iint_{\Pi^h} \big|\nabla\Phi_\e\big| dx =
\dfrac{2\pi}{|\ln \e|}\int_\e^{\sqrt \e} \bigg|\psi'\Big(\dfrac{\ln r}{\ln\e}\Big)\bigg| dr
\leqslant \dfrac{2\pi B_1 \sqrt \e}{|\ln\e|}
\end{equation}
и
\begin{equation}
 \label{eq-ph2}
\iint_{\Pi^h} \big|\nabla\Phi_\e\big|^2 dx =
\dfrac{2\pi}{|\ln \e|^2}\int_\e^{\sqrt \e} \bigg|\psi'\Big(\dfrac{\ln r}{\ln\e}\Big)\bigg|^2 \dfrac{dr}{r}
\leqslant \dfrac{2\pi B_1^2}{|\ln\e|^2} \int_\e^{\sqrt \e}  \dfrac{dr}{r}= \dfrac{\pi B_1^2}{|\ln\e|}.
\end{equation}

Now let $l\in\NN$. Let us pick the eigenfunctions $u^k_j$, $j=1,\dots,l$,
of the operator $H^0(k)$ associated respectively with the eigenvalues $E^0_j(k)$, i.e.
$H^0(k)u^k_j=E^0_j(k)u^k_j$, in such a way that the normalization conditions
\[
\langle u^k_j,u^k_s\rangle_{L^2(\Pi^h)}=\begin{cases}
1,& j=s,\\
0,& j\ne s.
\end{cases}
\]
hold. Note that the restrictions of the eigenfunctions $u^k_j$ onto the rectangles $\Pi^h_\pm$
are either zero functions or given by the explicit expressions \eqref{eq-ef}.
It follows that one can choose a constant $B_2>0$ such that
\[
\sup_{x\in \Pi^h_+\cup \Pi^h_-} \big|u^k_j(x)\big|+ \big|\nabla u^k_j(x)\big|\leqslant B_2, \quad
j=1,\dots,l, \quad k\in(-\pi,\pi],
\]
which can be used to prove that there exists a constant $B_3>0$ for which
\begin{equation}
    \label{eq-bnorm}
\Big|\langle \Phi_\e u^k_j,\Phi_\e u^k_s\rangle_{L^2(\Pi^h)}
-\langle u^k_j,u^k_s\rangle_{L^2(\Pi^h)}\Big|
\leqslant B_3 \e
\end{equation}
for all $j,s=1,\dots,l$ и всех $k\in(-\pi,\pi]$.
In particular, it follows from this last estimate that the
$l$ functions $\Phi_\e u^k_j$, $j=1,\dots, l$ are linearly independent for small $\e>0$.
Moreover, as noted above, each of these functions belongs to в $D^\e$.

Let us consider now an arbitrary subspace $L\in S_{l-1}$.
From the dimension considerations, the intersection $D^\e\cap L$ should contain at least one linear combination
$v$ of the form
\[
v=\sum_{j=1}^l b_j \Phi_\e u^k_j, \quad b=(b_j)\in\CC^l\setminus\{0\},
\]
and, due to \eqref{eq-bnorm}, for any vector $b$ there holds
\[
\|v\|^2_{L^2(\Pi^h)}\geqslant \|b\|^2_{\CC^l}(1- B_4\e)
\]
with some constant $B_4>0$ which is independent of $\e$ and $k$.
By elementary computations we obtain
\begin{multline*}
a^0_k(v,v)=\sum_{j,s=1}^l \overline{b_j}b_s
a^0_k(\Phi_\e u^k_j,\Phi_\e u^k_s)\\
=\sum_{j,s=1}^l \overline{b_j}b_s\bigg\{
\int_{\Pi^h} \Phi_\e^2(x) \Big[
\overline{\big(i\dfrac{\partial}{\partial x_1}-\dfrac{k}{2h}\big)u^k_j}
\big(i\dfrac{\partial}{\partial x_1}-\dfrac{k}{2h}\big)u^k_s
+\overline{\dfrac{\partial u^k_j}{\partial x_2}}\dfrac{\partial u^k_s}{\partial x_2}
\Big] dx
+ \int_{\Pi^h} |\nabla \Phi_\e|^2 \overline{u^k_j}u^k_s\, dx\\
+\int_{\Pi^h} \bigg[
\Phi_\e\dfrac{\partial\Phi_\e}{\partial x_1}
\Big(
\overline{\big(i\dfrac{\partial u^k_j}{\partial x_1}-\dfrac{ku^k_j}{2h}\big)}u^k_s
+\overline{u^k_j} \big(i\dfrac{\partial u^k_s}{\partial x_1}-\dfrac{k u^k_s}{2h}\big)
\Big)+
\Phi_\e \dfrac{\partial\Phi_\e}{\partial x_2}\Big(
\overline{\dfrac{\partial u^k_j}{\partial x_2}}u^k_s
+
\overline{u_j^k}\dfrac{\partial u^k_s}{\partial x_2}
\Big)
\bigg]\,dx
\end{multline*}

It follows, by using \eqref{eq-ph1} and \eqref{eq-ph2} that for all
$j,s=1,\dots, l$ and $k\in(-\pi,\pi]$ there holds
\begin{gather*}
\Big|\int\limits_{\Pi^h} |\nabla \Phi_\e|^2 \overline{u^k_j}u^k_s\, dx\Big| \leqslant \dfrac{B_5}{|\ln\e|},\\
\bigg|\int\limits_{\Pi^h} \bigg[
\Phi_\e\dfrac{\partial\Phi_\e}{\partial x_1}
\Big(
\overline{\big(i\dfrac{\partial u^k_j}{\partial x_1}-\dfrac{k u^k_j}{2h}\big)}u^k_s
+\overline{u^k_j} \big(i\dfrac{\partial u^k_s}{\partial x_1}-\dfrac{k u^k_s}{2h}\big)
\Big)+
\Phi_\e \dfrac{\partial\Phi_\e}{\partial x_2}\Big(
\overline{\dfrac{\partial u^k_j}{\partial x_2}}u^k_s
+
\overline{u_j^k}\dfrac{\partial u^k_s}{\partial x_2}
\Big)
\bigg]\,dx
\bigg|
\leqslant \dfrac{B_6\sqrt\e}{|\ln\e|},\\
\bigg|\int\limits_{\Pi^h} \Phi_\e^2(x) \Big[
\overline{\big(i\dfrac{\partial}{\partial x_1}-\dfrac{k}{2h}\big)u^k_j}
\big(i\dfrac{\partial}{\partial x_1}-\dfrac{k}{2h}\big)u^k_s
+\overline{\dfrac{\partial u^k_j}{\partial x_2}}\dfrac{\partial u^k_s}{\partial x_2}
\Big] dx-a^0_k(u^k_j,u^k_s)\bigg|\leqslant B_7\e
\end{gather*}
with some positive constants $B_5,B_6,B_7$ independent of $k$ and $\e$.

By noting that $a^0_k(u^k_j,u^k_s)=E^0_j(k)\langle u^k_j,u^k_s\rangle_{L^2(\Pi^h)}$
and combining all the above estimates we arrive at
\[
a^0_k(v,v)\leqslant \sum_{j=1}^l E^0_j(k) |b_j|^2 + \dfrac{B_8\|b\|^2_{\CC^l}}{|\ln\e|},
\]
where the constant $B_8>0$ does not depend on $b$, $\e$ and $k$.
Now, due to $E^0_l(k)=\max_{j\in\{1,\dots,l\}}E^0_j(k)$, we have
\[
\min_{\substack{u\in D^\e\cap L\\v\ne 0}} \dfrac{a^0_k(u,u)}{\|u\|^2_{L^2(\Pi^h)}}
\leqslant
\dfrac{a^0_k(v,v)}{\|v\|^2_{L^2(\Pi^h)}}\leqslant
\dfrac{\sum_{j=1}^l E^0_j(k) |b_j|^2 + \dfrac{B_8\|b\|^2_{\CC^l}}{|\ln\e|}}{(1-B_4\e)\|b\|^2_{\CC^l}}
\leqslant E^0_l(k) + \dfrac{C_l}{|\ln\e|}
\]
with some constant  $C_l$ independent of $b,\e,k$.
As the subspace $L\in S_{l-1}$ was arbitrary, the right-hand side
of \eqref{eq-ele} follows from the max-min principle \eqref{eq-minmax}.

\section{Analysis near the crossings of band functions}\label{sec3}

In this section we proceed with a more detailed analysis of the eigenvalues $E^\e_l(k)$
near $k=k_0$ and conclude the proof of Theorem~\ref{th1.1}.

It follows directly from the definition of the functions $E_{m,0}^\pm$ that if the condition (\ref{1.4}) holds
for some $m$, $n$, $k_0$, then it is also valid for $-m$, $-n$, $-k_0$,
and the value $E_0$ remains unchanged.
Furthermore, there are no other points at which (\ref{1.4a}) holds with the same value $E_0$.
Therefore, it follows from Lemma~\ref{prop1} that
for any $A>0$ and sufficiently small $\e$ there exists a constant $t_0=t_0(A)>0$ for which the inequality
\begin{equation}\label{3.39}
\dist(\spec(H_\e(k)),E_0)\geqslant \frac{A}{|\ln\e|}\quad \text{for}\quad |k-k_0|\geqslant \frac{t_0}{|\ln\e|}
\text{ and }  |k+k_0|\geqslant \frac{t_0}{|\ln\e|}
\end{equation}
holds.
Therefore, for those $k$ which are sufficiently far away from $\pm k_0$, the spectrum of $H^\e(k)$ is separated from the point $E_0$
by a distance of at least $\Odr(|\ln\e|^{-1})$. For the proof of theorem it is now necessary to investigate
the behavior of the spectrum of the operator $H^\e(k)$ as $|k-k_0|\leqslant t_0 |\ln\e|^{-1}$ and $|k+k_0|\leqslant t_0|\ln\e|^{-1}$.
Such a study will be done in a neighborhood of the point $k_0$ only. For the other case,
the results can be easily obtained using the obvious fact that
the replacement $k\mapsto -k$ in the operator $H^\e(k)$ does not change the eigenvalues,
and the associated eigenfunctions are complex conjugate to the original ones.

We will consider the values $k$ having the form
\begin{equation*}
k_\e=k_0+\frac{t}{\ln\e},
\end{equation*}
where $t$ is a real parameter running through a segment $[-t_0,t_0]$. By Lemma~\ref{prop1}, the operator $H^\e(k_\e)$ has exactly
two eigenvalues converging to $E_0$ for $\e\to+0$. Denote them by $E^{(j)}_\e(t)$, $j=1,2$. We will assume the ordering
$E^{(1)}_\e(t)\leqslant E^{(2)}_\e(t)$.

To construct the asymptotics of $E^{(j)}_\e(t)$ for  $\e\to0$ we will use the methods of matched asymptotic expansions \cite{Il}.
This will be done in two steps: first we construct the formal asymptotic expansions,
and then they will be rigorously justified.

\subsection{Constructing formal asymptotics}

We will construct the asymptotics of the above eigenvalues in the following form:
\begin{equation}\label{3.5}
E^{(j)}_\e(t)=E_0+\frac{\mu^{(j)}(t)}{\ln\e}+\ldots,
\end{equation}
where $\mu^{(j)}(t)$ are some constants and should be determined.

Outside a small neighborhood of the window $\varpi^\e_0$
the asymptotics of the eigenfunctions $\psi_\e^{(j)}$ of $H^\e(k_\e)$ associated to the eigenvalues $E^{(j)}_\e$
will be sought in the form
\begin{equation}
\psi_\e^{(j)}(x)=\E^{\iu \frac{2\pi m}{2h}x_1} \left( \phi_0^{(j)}(x)+\frac{\phi_1^{(j)}(x)}{\ln \e}+\ldots
\right),
\label{3.6}
\end{equation}
with
\begin{equation}   \label{3.7}
\phi_0^{(j)}=
\begin{cases}
a_+^{(j)}\phi_0^+(x), & x_2>0,
\\
a_- ^{(j)}\phi_0^-(x), & x_2<0,
\end{cases}
\qquad
\begin{aligned}
\phi_0^+(x)&=\E^{\iu\frac{\pi (n-m)}{h}\,x_1} \cos \frac{x_2}{2},\\
\phi_0^-(x)&=\cos\frac{\pi x_2}{2d},
\end{aligned}
\end{equation}
where $a_\pm^{(j)}$ are constants which cannot vanish simultaneously.
Everywhere below we call the representation (\ref{3.6})
the outer expansion.

Let us substitute the outer expansion
and (\ref{3.5}) into the eigenvalue equation
\begin{equation}\label{3.8}
H^\e(k)\psi^{(j)}_\e=E_\e^{(j)}\psi^{(j)}_\e,
\end{equation}
and then group the coefficients  at the like powers of $\dfrac{1}{\ln \e}$.
Considering the first two terms of the series obtained gives the following
problems for the coefficients $\phi_1^{(j)}$:
\begin{equation}\label{3.9}
\begin{aligned}
&\left(\Big(\iu \frac{\p}{\p x_1} -\frac{k_0+2\pi m}{2h}\Big)^2 -\frac{\p^2}{\p x_2^2}-E_0\right)\phi_1^{(j)}
\\
&\hphantom{\left(\iu \frac{\p}{\p x_1} -\frac{k+2\pi m}{h}\right)^2}
=\left(\mu^{(j)} + \frac{t}{h} \Big(\iu \frac{\p}{\p x_1} - \frac{k_0+2\pi m}{2h}\Big)
\right)\phi_0^{(j)} \quad\text{in}\quad \Pi^h,
\\
&\frac{\p\phi_1^{(j)}}{\p x_2}=0\quad\text{on}\quad \G_\mathrm{N}^{h,0},\qquad \phi_1^{(j)}=0\quad \text{on} \quad \G^h_\mathrm{D},
\end{aligned}
\end{equation}
and the periodic boundary conditions (\ref{1.10c}) for the functions  $\phi_1^{(j)}$
are imposed at the side boundary of the domain $\Pi^h$.
The conditions obtained split into two independent problems on $\Pi^h_+$ and $\Pi^h_-$.
Their solvability in $\H^1(\Pi^h_\pm)$ is equivalent to the orthogonality of the right-hand sides to the the functions
$\phi_0^\pm$ in $L_2(\Pi^h_\pm)$. It is easy to see that the solvability conditions will be satisfied only if $a_+^{(j)}=0$ or $a_-^{(j)}=0$.
On the other hand, in that case one cannot match the asymptotic expansions near the window $\omega^\e_0$.
Such difficulties are quite standard for problems involving a change of boundary conditions on a small part of the boundary,
see e.g. \cite[Chs. 2 and 6]{GN} or
\cite[Ch. I\!I\!I, \S 1]{Il}, and they are usually overcome by considering the problem (\ref{3.9})
in a larger sense.  Assume that $\phi_1^{(j)}$ has a logarithmic singularity at the contraction point of the window $\omega^\e_0$:
\begin{equation}\label{3.10}
\phi_1^{(j)}(x)=b_\pm^{(j)} \ln |x|+ \widetilde{b}_\pm^{(j)} + \Odr(|x|), \quad x\to0,\quad \pm x_2>0,
\end{equation}
where $b_\pm^{(j)}$, $\widetilde{b}_\pm^{(j)}$ are some constants,
and we assume that that the function $\phi_1^{(j)}$ is sufficiently regular
outside the point $x=0$ and satisfies all the boundary conditions mentioned above.

Let us multiply the equations (\ref{3.9}) by $\overline{\phi_0^\pm}$ and then integrate
by parts in $\Pi^h_\pm$ taking into account the asymptotics (\ref{3.10}).
This gives necessary solvability conditions for the problems (\ref{3.9}):
\begin{equation}\label{3.11}
\begin{aligned}
\pi b_+^{(j)}&=a_+^{(j)}\left( \pi h \mu^{(j)} -\frac{\pi t}{2h} (2\pi n+k_0)\right),
\\
\pi b_-^{(j)}&= a_-^{(j)}\left(  h \mu^{(j)} d -\frac{t d}{2h} (2\pi m+k_0)\right).
\end{aligned}
\end{equation}
Using the standard arguments, see e.g. \cite[Ch. I\!I\!I, \S 2, Th.~2.2]{Il}
it is possible to show that these conditions are also sufficient.
Note that there are no restrictions on the coefficients $\widetilde{b}_\pm^{(j)}$,
and one can choose them freely by adding to $\phi_1^{(j)}$ the functions of the form $c_\pm \phi_0^\pm$.

The first terms of the expansions (\ref{3.6}) constructed above
may serve as an approximation of the true eigenfunctions
outside a certain neighborhood of the window $\varpi^\e_0$ only,
while the boundary conditions near the window are not satisfied.
To construct the asymptotics near the window let us introduce the stretched coordinates
$\xi=(\xi_1,\xi_2)=x\e^{-1}$ and look for the asymptotics of the eigenfunctions
in the form of an inner expansion:
\begin{equation}\label{3.12}
\psi_\e^{(j)}(x)=\E^{\iu \frac{2\pi m}{2h}x_1} \left(
w_0^{(j)}(\xi)+\frac{1}{\ln\e} w_1^{(j)}(\xi)+\ldots
\right)
\end{equation}
Let us substitute this ansatz and (\ref{3.5}) into the equation (\ref{3.8})
and group the coefficients in front of $\e^{-2}$ and $\e^{-2}\ln^{-1}\e$ in this formally obtained equality.
This results in the following conditions on the functions $w_p^{(j)}$, $p=1,2$:
\begin{equation}\label{3.13}
\begin{aligned}
&\D_\xi w_p^{(j)}=0\quad\text{in}\quad \mathds{R}^2\setminus\overline{\G}_*, \quad  \frac{\p w_p^{(j)}}{\p \xi_2}=0\quad\text{on}\quad \overline{\G}_*, \quad \G_*:=O\xi_1\setminus[-1,1].
\end{aligned}
\end{equation}
To determine uniquely the functions  $w_p^{(j)}$ we need to prescribe their behavior at infinity.
This can be found from the matching of the inner expansion  (\ref{3.12}) with the outer expansion (\ref{3.6}).
Let us represent the functions $\phi_0^{(j)}$ and $\phi_1^{(j)}$
in the variables $\xi$ in a neighborhood of $x=0$:
\begin{align*}
\phi_0^{(j)}(x)&+\frac{1}{\ln\e} \phi_1^{(j)}(x)=a_\pm^{(j)}+b_\pm^{(j)}
\\
&+ \frac{1}{\ln\e} (b_\pm^{(j)}\ln|\xi|+\widetilde{b}_\pm^{(j)})+\Odr\big(\e(|\xi|+1)\big),\quad x\to0,\quad \pm x_2>0.
\end{align*}
By matching these relations with (\ref{3.12}),
we obtain
\begin{align}
&w_0^{(j)}(\xi)=a_\pm^{(j)}+b_\pm^{(j)}+o(1),\quad |\xi|\to\infty,\quad \pm\xi_2>0,
\label{3.15}
\\
& w_1^{(j)}(\xi)= b_\pm^{(j)} \ln|\xi|+\widetilde{b}_\pm^{(j)}, \hphantom{1)} \quad |\xi|\to\infty,\quad \pm\xi_2>0.\label{3.16}
\end{align}
Therefore, one needs to find the solutions of the problems (\ref{3.13})
behaving at infinity in the prescribed way.

The solvability of the problems (\ref{3.13}) is rather easy to study with the help
of the methods of complex analysis. Namely, the domain $\mathds{R}^2\setminus\overline{\G}_*$ can be conformally mapped
onto the disk with the help of an explicitly given conformal map. The problem in the disk
can be then studied in a standard way, one just needs to take into account the behavior at the image of infinity.
Such an approach gives the following results on the solvability of the problems (\ref{3.13}), (\ref{3.15}), (\ref{3.16}).

The problem (\ref{3.13}) with the infinity behavior (\ref{3.15}) has the unique solution which is the constant one:
$w_0^{(j)}(\xi)\equiv c_0$.
This is only possible if the solvability conditions
\begin{equation}\label{3.18}
a_+^{(j)}+b_+^{(j)}=a_-^{(j)}+b_-^{(j)}=c_0
\end{equation}
hold. The problem  (\ref{3.13}) with the condition (\ref{3.16}) can also be solved uniquely (up to a constant)
under the compatibility conditions
\begin{equation}\label{3.19}
b_+^{(j)}=-b_-^{(j)}.
\end{equation}
Its solution can be constructed explicitly:
\begin{equation}\label{3.20}
\begin{aligned}
&w_1^{(j)}(\xi)=b_+^{(j)} X(\xi)+c_1,\quad c_1=\mathrm{const},
\\
&X(\xi):=\RE\ln (z+\sqrt{z^2-1}),\quad z=\xi_1+\iu \xi_2.
\end{aligned}
\end{equation}
Here the branches of the square root and the logarithm are chosen by the conditions $\sqrt{1}=1$, $\ln 1=0$.

Let us emphasize some properties of the function $X$, which can be easily checked using its explicit definition. The function $X$
is infinitely smooth in $\mathds{R}^2\setminus\overline{\G}_*$ and is contiuous up to the boundary if this domain.
It also has one-side derivatives on the set $\G_*$, which should be understood as a two-side cut.
The function $X$ is odd with respect to $\xi_2$ and has the following asymptotic behavior at infinity:
$X(\xi)=\pm (\ln|\xi|+\ln 2)+\Odr(|\xi|^{-2})$, $|\xi|\to\infty$, $\pm\xi_2>0$.
This last equality leads to the condition (\ref{3.19}). Due to the asymptotics (\ref{3.16}), this equality also
allows one to determine the constants $\widetilde{b}_\pm^{(j)}$:
$\widetilde{b}_\pm^{(j)}=b_\pm^{(j)}\ln 2$.

By elementary computations one can easily check that the system (\ref{3.18}), (\ref{3.19}) is equivalent to the following one:
\begin{equation*}
\left\{
\begin{aligned}
-2 b_+^{(j)} - a_+^{(j)} + a_-^{(j)}&=0,
\\
\frac{\pi}{d} a_+^{(j)} - \frac{\pi}{d} (a_-^{(j)}+ 2 b_-^{(j)})&=0.
\end{aligned}
\right.
\end{equation*}
Let us express the numbers  $b_\pm^{(j)}$ in terms of $\mu^{(j)}$ using the formulas (\ref{3.11})
and substitute the expressions obtained into the last system of equations:
\begin{gather}
(M-2h\mu^{(j)}) a^{(j)} =0, \nonumber
\\
M:=
\begin{pmatrix}
& t\b_1-1 & 1
\\
& \z & -\z-t\b_2
\end{pmatrix},
\quad a^{(j)}:=
\begin{pmatrix}
a_+^{(j)}
\\
a_-^{(j)}
\end{pmatrix},\label{3.24a}
\end{gather}
where
\begin{equation}\label{3.24b}
\b_1:= \frac{2\pi n+ k_0}{h},\quad \b_2:=-\frac{2\pi m+k_0}{h}, \quad
\zeta=\dfrac{\pi}{d}.
\end{equation}
As we assume that the vectors $a^{(j)}$ are non-zero, the quantities $2h \mu^{(j)}$ must be the eigenvalues
of the matrix $M$, therefore,
\begin{gather}
\mu^{(j)}(t)=\frac{f_j(t,\z)}{4h},\label{3.25}
\\
\begin{aligned}
&f_1(t,\z):=t(\b_1-\b_2)-(\z+1)+ \sqrt{\big(
t(\b_1+\b_2)+\z-1\big)^2+4\z},
\\
&f_2(t,\z):=t(\b_1-\b_2)-(\z+1)- \sqrt{\big(
t(\b_1+\b_2)+\z-1\big)^2+4\z}.
\end{aligned}\nonumber
\end{gather}
The respective eigenvectors can be chosed by:
\begin{equation}\label{3.26}
a^{(j)}=
\begin{pmatrix}
1
\\
2h\mu^{(j)}+1-t\b_1
\end{pmatrix}, \quad j=1,2.
\end{equation}
The expressions obtained for $a_\pm^{(j)}$
completely determine the inner and the outer expansions for the eigenfunctions $\psi_\e^{(j)}$.
The constant $c_1$ in (\ref{3.20}) will be set equal to zero. Therefore, the formal construction
of the asymptotics is completed.

\subsection{Justifying the formal asymptotics}

The first step of the justfication consists in a proof of the fact that
the asymptotics (\ref{3.5}), (\ref{3.6}) deliver formal asymptotic solutions of the equation (\ref{3.8}), i.e.
satisfy it with some small error. In our case one needs to obtain an error of order $\Odr\Big(\dfrac{1}{|\ln\e|^2}\Big)$.
For that, it is necessary to construct some additional terms in the outer and the inner expansions of the eigenfunctions $\psi_\e^{(j)}$
and some additional summands in the asymptotics (\ref{3.5}). This is a standard situation in problems with singular perturbations.
In our case, such a construction does not involve any conceptual difficulties but is rather technical.
Furthermore, it consists essentially in a repeating of the similar constructions of \cite{GN}, \cite[Ch. I\!I\!I, \S 1]{Il}.
Hence, we will not enter into details of the construction of additional terms, but prefer
to formulate directly the final result. Let us denote
\begin{align*}
\widehat{\psi}_\e^{(i)}(x)=&\E^{\iu \frac{2\pi m}{h} x_1}
\left[\left(\phi_0^{(j)}(x)+ \frac{1}{\ln \e} \phi_1^{(j)}(x)
\right)\chi \left(\frac{|x|}{\e^\a}\right)\right.
\\
&\hphantom{\E^{\iu \frac{k+2\pi m}{h} x_1}[
}\left.+ \left(w_0^{(j)}(\xi)+\frac{1}{\ln\e}
\right) \left(1-\chi \left(\frac{|x|}{\e^\a}\right)\right)
\right]
\\
\widehat{E}_\e^{(j)}(t):=&E_0 + \frac{\mu^{(j)}(t)}{\ln\e},
\end{align*}
whereе $\chi=\chi(s)$ is a smooth function which a equal to zero for $s<1$ and equal to one for $s>2$,
and $\a$ is some constant.

\begin{lemma}\label{lm3.1}
The function $\widehat{\psi}_\e^{(j)}$ belongs to $C^\infty(\Pi^h)\cap \H^1(\Pi^h)$. One has the convergence
\begin{equation}\label{3.28}
\widehat{\psi}_\e^{(j)}\to \E^{\iu \frac{2\pi m}{2h}x_1}\phi_0^{(j)} \quad \text{in}\quad L_2(\Pi^h),\qquad \widehat{E}_\e^{(j)}\to E_0.
\end{equation}
There exists a number $\a\in(0,1)$, functions $\psi_*^{(j)}(x,\e)\in C^\infty(\Pi^h) \cap \H^1(\Pi^h)$, and numbers $E_*^{(j)}(\e)$
such that the functions
\begin{equation*}
\widehat{\psi}_{\e,*}^{(j)}(x):=\widehat{\psi}_\e^{(j)}(x)+\psi_*^{(j)}(x,\e),
\quad \widehat{E}_{\e,*}:=E_\e^{(j)}+ E_*^{(j)}(\e)
\end{equation*}
satisfy the equation
\begin{equation}\label{3.28a}
(H^\e(k)-\widehat{E}_{\e,*}^{(j)})\widehat{\psi}_{\e,*}^{(j)}=g_\e^{(j)},
\end{equation}
one has, uniformly in $\e$ and $t$, the estimates
\begin{align}
&\|\psi_*^{(j)}(\cdot,\e)\|_{\H^1(\Pi^h)}=\Odr(\ln^{-2}\e),\quad
|E_*^{(j)}(\e)|=\Odr(\ln^{-2}\e),\label{3.29}
\\
&\|g_\e^{(j)}\|_{L_2(\Pi^h)}=\Odr(\ln^{-2}\e).\label{3.30}
\end{align}
\end{lemma}

Now let us consider the remainders of the asymptotics. By (\ref{3.28a}) and \cite[Ch. V, \S 3.5, Eq. (3.21)]{Ka}
there holds
\begin{equation}\label{3.33}
\widehat{\psi}^{(j)}_{\e,*}= \big(H^\e(k_\e)-\widehat{E}_*^{(j)}(\e)\big)^{-1} g_\e^{(j)}=\sum\limits_{p=1}^{2} \frac{(g_\e^{(j)},\psi_\e^{(p)})_{L_2(\Pi^h)}}{E_\e^{(p)}-\widehat{E}_\e^{(j)}}
\psi_\e^{(p)} + R_\e^{(j)} g_\e^{(j)},
\end{equation}
where, $E_\e^{(p)}$ are the eigenvalues of the operator $H^\e(k)$ converging to $E_0$ as $\e\to+0$, $\psi_\e^{(p)}$ are
the respective eigenfunctions which are supposed to be orthonormal in $L_2(\Pi^h)$, and $R_\e^{(j)}$ is the reduced resolvent
considered as an operator in $L_2(\Pi^h)$ and satisfying the estimate
\begin{equation*}
\|R_\e^{(j)}\|=\frac{1}{\dist\big(\widehat{E}_\e^{(j)},\spec(H^\e(k_\e))\setminus
\{E_\e^{(1)},E_\e^{(2)}\}\big)}\leqslant C,
\end{equation*}
where the constant $C$ does not depend on $\e$ and $t$. Furthermore, the operator
$R_\e^{(j)}$ acts in the orthogonal complement of $\{\psi_\e^{(1)}, \psi_\e^{(2)}\}$ in $L_2(\Pi^h)$.

Let us denote
\begin{equation*}
A_{jp}(\e):=\frac{(g_\e^{(j)},\psi_\e^{(p)})_{L_2(\Pi^h)}} {E_\e^{(p)}-E_\e^{(j)}}, \quad A_j(\e):=(A_{j1}(\e),A_{j2}(\e)).
\end{equation*}
It follows from the properties of the operators $R_\e^{(j)}$ described above and from (\ref{3.33}), (\ref{3.30}), (\ref{3.28}), (\ref{3.26}), (\ref{3.24b}), (\ref{3.24a}), (\ref{3.7}) that
\begin{align*}
&\|\widehat{\psi}_{\e,*}^{(j)}\|_{L_2(\Pi^h)}^2= \sum\limits_{p=1}^{2} |A_{jp}|^2+ \|R_\e^{(j)} g_\e^{(j)}\|_{L_2(\Pi^h)}^2,\quad j=1,2,
\\
& (\widehat{\psi}_\e^{(1)}, \widehat{\psi}_\e^{(2)})_{L_2(\Pi^h)} = \sum\limits_{p=1}^{2}
A_{1p} A_{2p} + (R_\e^{(1)} g_\e^{(1)},R_\e^{(2)} g_\e^{(2)})_{L_2(\Pi^h)},
\end{align*}
and
\begin{equation*}\begin{aligned}
&\|A_j\|_{\mathds{R}^2}=\|\phi_0^{(j)}\|_{L_2(\Pi^h)}+o(1), \quad
\e\to+0,
\\
&(A_1,A_2)_{\mathds{R}^2}=(\phi_0^{(1)},\phi_0^{(2)})_{L_2(\Pi^h)}
+o(1)=o(1),\quad
\e\to+0,
\end{aligned}
\end{equation*}
uniformly in $t\in[-t_0,t_0]$. At cost of redenoting $\psi_\e^{(1)}$ by $\psi_\e^{(2)}$ and vice versa,
due to the above relations, we conclude that one has the estimates
\begin{equation*}
|A_{11}|\geqslant C>0,\quad |A_{22}|\geqslant C>0,
\end{equation*}
which are uniform in $\e$ and $t\in[-t_0,t_0]$.
Therefore, by the definition of the functions $A_{jp}$ and by the estimate (\ref{3.30}), we obtain
the estimates of the remainder
\begin{equation*}
E_\e^{(j)}(t)-\widehat{E}_{\e,*}^{(j)}(t)=\Odr(\ln^{-2}\e)
\end{equation*}
uniformly in $t\in[-t_0,t_0]$.
Now, taking into account (\ref{3.29}), we obtain the final form of the asymptotics
for the eigenvalues  $E_\e^{(j)}$ uniformly in $t\in[-t_0,t_0]$:
\begin{equation}\label{3.36}
E_\e^{(j)}(t)=E_0+ \frac{\mu^{(j)}(t)}{\ln\e}+\Odr(\ln^{-2}\e),\quad j=1,2.
\end{equation}

At the end of the section let us remark that, using the methods of the works cited above,
one can construct an asymptotics of the eigenvalues $E_\e^{(j)}$ up to a certain power of $\e$
by summing up the logarithmic part. On the other hand, such a precision is redundant for our objectives.

\subsection{Studying the asymptotic formulas}

Let us proceed with the qualitative analysis of the band functions near $k=k_0$
using the asymptotics obtained. To start, let us determine the extrema of the functions $f_j$ in (\ref{3.25}).

It follows from the inequality (\ref{1.4a}) that the numbers
$\b_1$ and $\b_2$ in (\ref{3.24b}) have the same sign. Hence
\begin{equation*}
|\b_1-\b_2|=\big||\b_1|-|\b_2|\big|,\quad |\b_1+\b_2|=|\b_1|+|\b_2|,
\end{equation*}
which implies $|\b_1-\b_2|<|\b_1+\b_2|$.
Therefore,
\begin{equation*}
\lim\limits_{t\to\pm\infty} f_1(t,\z)=+\infty,\quad \lim\limits_{t\to\pm\infty} f_2(t,\z)=-\infty.
\end{equation*}
Using the above and the continuity of the functions $f_1$ and $f_2$ with respect to $t$ one sees
that they attain their respective global minimum and maximum at some finite points.
Calculating these points in the standard way we obtain
\begin{gather*}
\min\limits_{t\in \mathds{R}} f_1(t,\z)=f_1(t_{\mathrm{min}},\z),\quad \max\limits_{t\in\mathds{R}} f_2(t,\z)=f_2(t_{\mathrm{max}},\z),\label{3.40}
\\
\begin{aligned}
&t_{\mathrm{min}}=-\frac{2\b}{\b_1+\b_2}
\sqrt{\frac{\z}{1-\b^2}}+\frac{1-\z}{\b_1+\b_2},
\\
&
t_{\mathrm{max}}=\frac{2\b}{\b_1+\b_2}\sqrt{\frac{\z}
{1-\b^2}} +\frac{1-\z}{\b_1+\b_2},
\end{aligned}
\end{gather*}
where
\begin{equation*}
\b:= \frac{\b_1-\b_2}{\b_1+\b_2}=\frac{\pi(m+n)+k_0}{\pi(n-m)}.
\end{equation*}
The values of the global minimum of the function $f_1$ and of the global maximum of the function $f_2$
are given by the formula
\begin{equation}\label{3.42}
\begin{aligned}
&f_1(t_{\mathrm{min}},\z)=2\sqrt{\z(1-\b^2)}-\b(\z-1)-\z-1=\tau_+,
\\
&f_2(t_{\mathrm{max}},\z)=-2\sqrt{\z(1-\b^2)}-\b(\z-1)-\z-1=\tau_-.
\end{aligned}
\end{equation}
The inequality $\ln\e<0$ together with the preceding relations and the equalities (\ref{3.25}), (\ref{3.36}), (\ref{3.40}) give
\begin{equation*}
\begin{aligned}
\max\limits_{t\in[-t_0,t_0]} E_\e^{(1)}=&E_0- \frac{1}{4h|\ln\e|}\min\limits_{t\in[-t_0,t_0]} f_1(t,\z)+\Odr\Big(\dfrac{1}{|\ln\e|^2}\Big)
\\
=&E_0-\frac{\tau_+}{4h|\ln\e|}+\Odr\Big(\dfrac{1}{|\ln\e|^2}\Big),
\\
\min\limits_{t\in[-t_0,t_0]} E_\e^{(2)}=&E_0- \frac{1}{4h|\ln\e|}\max\limits_{t\in[-t_0,t_0]} f_2(t,\z)+\Odr\Big(\dfrac{1}{|\ln\e|^2}\Big)
\\
=&E_0-\frac{\tau_-}{4h|\ln\e|}+\Odr\Big(\dfrac{1}{|\ln\e|^2}\Big),
\end{aligned}
\end{equation*}
where the parameter $t_0$ can be chosen sufficiently large to have $t_{\mathrm{min}}\in[-t_0,t_0]$,  $t_{\mathrm{max}}\in[-t_0,t_0]$.
It follows from the last equality and from (\ref{3.42}) that, for $\e$ sufficiently small,
\begin{equation*}
\begin{aligned}
\min\limits_{t\in[-t_0,t_0]} E_\e^{(2)}(t)- \max\limits_{t\in[-t_0,t_0]} E_\e^{(1)}(t)=& \frac{\tau_+- \tau_-}{|\ln\e|}+\Odr\Big(\dfrac{1}{|\ln\e|^2}\Big)
\\
=&\frac{\sqrt{z(1-\b^2)}}{h|\ln\e|} +\Odr\Big(\dfrac{1}{|\ln\e|^2}\Big)>0.
\end{aligned}
\end{equation*}
Taking into account the inequality (\ref{3.39}) and choosing the parameter $A$ sufficiently large we conclude
that there is a gap in the spectrum of the operator $H^\e$,
\begin{equation*}
\left( E_0-\frac{\tau_+}{4h|\ln\e|}+\Odr\Big(\dfrac{1}{|\ln\e|^2}\Big), \  E_0- \frac{\tau_-}{4h|\ln\e|}+\Odr\Big(\dfrac{1}{|\ln\e|^2}\Big)\right),
\end{equation*}
which proves (\ref{1.3}) and (\ref{1.6}). It remains to show the validity of the asymptotics (\ref{1.7}).

As the function $f_1(t,\z)$ attains its global maximum at the point $t_{\mathrm{min}}$ only, for all $t$ at a finite distance from
$t_{\mathrm{min}}$ there holds $f_1(t,\z)>f_1(t_{\mathrm{min}},\z)$. Therefore, the point of the maximum of the function
$E_1^{(\e)}$ should tend to $t_{\mathrm{min}}$ as $\e\to+0$. As
\begin{equation*}
f_1(t,\z)=f_1(t_{\mathrm{min}},\z)+\Odr\big((t-t_{\mathrm{min}})^2
\big),\quad t\to t_{\mathrm{min}},
\end{equation*}
the equality
\begin{equation*}
E_\e^{(1)}(t)=E_0-\frac{\tau_+}{4h|\ln\e|} +\Odr\Big(\dfrac{1}{|\ln\e|^2}\Big)
\end{equation*}
can only be satisfied for $t$ situated at the distance of at most
$\Odr\Big(|\ln\e|^{-1/2}\Big)$ from $t_{\mathrm{min}}$.
This gives the asymptotics (\ref{1.7}) for $k_+^+(\e)$.
The asymptotics for $k^+_-(\e)$ is proved in the same way.

\section*{Acknowledgments}

The present work was done during the visit of D. I. Borisov to the University Paris-Sud in summer 2011
within the project ``Th\'eorie des op\'erateurs et probl\`emes au bord en m\'ecanique quantique''.
D. I. Borisov thanks  for the warm welcome and the hospitality shown to him.

The authors thank the referee, who remained anonymous, for a series of useful remarks
that helped to improve significantly the section  \ref{sec2}.

\end{document}